\documentclass[12pt]{article}

\usepackage{float}
\usepackage{geometry}
\usepackage{verbatim}
\usepackage{graphicx}
\usepackage{enumerate, enumitem}
\usepackage{subcaption}
\usepackage{adjustbox}
\usepackage{booktabs}
\usepackage[
colorlinks=true,
citecolor=blue,
linkcolor=magenta,
urlcolor=cyan]{hyperref}
\usepackage{url}
\usepackage{xcolor}
\usepackage{amssymb,amsthm,amsmath,amsfonts,amsthm,mathtools}
\newcommand{\eg}{{\it e.g.}}
\newcommand{\ie}{{\it i.e.}}

\newcommand{\BEQ}{\begin{equation}}
\newcommand{\EEQ}{\end{equation}}
\newcommand{\BEAS}{\begin{eqnarray*}}
\newcommand{\EEAS}{\end{eqnarray*}}
\newcommand{\ones}{\mathbf 1}
\newcommand{\reals}{{\mbox{\bf R}}}

\newcommand{\var}{\mathop{\bf var}}
\newcommand{\avg}{\mathop{\bf avg}}

\newcommand{\argmax}{\mathop{\rm argmax}}

\newcounter{algorithmctr}
\renewcommand{\thealgorithmctr}{\arabic{algorithmctr}}
\newenvironment{algdesc}%
   {\mbox{}\\*[\parskip]\begin{minipage}{\linewidth}%
       \refstepcounter{algorithmctr}\begin{list}{}{%
       \setlength{\rightmargin}{0\linewidth}%
       \setlength{\leftmargin}{.05\linewidth}}%
       \rmfamily\small
       \item[]{\setlength{\parskip}{0ex}\hrulefill\par%
        \nopagebreak{\bfseries\textsf{Algorithm \thealgorithmctr~}}}}%
   {{\setlength{\parskip}{-1ex}\nopagebreak\par\hrulefill\\*[2ex]\par}%
   \end{list}\end{minipage}}

\usepackage{textcomp}
\usepackage{listings}
\definecolor{seagreen}{rgb}{0.18, 0.55, 0.34}
\definecolor{mediumviolet-red}{rgb}{0.78, 0.08, 0.52}
\definecolor{khaki}{rgb}{0.94, 0.9, 0.55}

\lstdefinelanguage{mypython}
{
	keywords=[1]{from, import, assert, not, print},
	keywordstyle=[1]{\color{mediumviolet-red}},
	keywords=[2]{resalloc, fungible, AllocationProblem, utilities, torch},
	keywordstyle=[2]{\color{seagreen}},
	numbers=none,
	upquote=true,
	showstringspaces=false,
	basicstyle=\ttfamily,
	columns=fullflexible,
	keepspaces=true,
	emph={True,False,as,def,return,float},
	emphstyle={\color{seagreen}},
	frame=trBL,
	belowskip=1em,
	aboveskip=1em,
	captionpos=b
}

\usepackage[
natbib=true,
citestyle=alphabetic,
style=alphabetic,
maxnames=3,
minnames=3,
maxbibnames=99,
backend=bibtex,
urldate=iso8601,
firstinits=true,
isbn=false,
doi=false,
date=year,
sorting=nty
]{biblatex}
\DeclareFieldFormat{eprint:arXiv}{%
{\hypersetup{hidelinks}\href{https://www.arXiv.org/abs/#1}{arXiv:#1}}%
}
\DeclareFieldFormat[article,inbook,incollection,inproceedings,patent,thesis,
unpublished]{citetitle}{#1}
\DeclareFieldFormat[article,inbook,incollection,inproceedings,patent,thesis,
unpublished]{title}{#1}
\renewbibmacro{in:}{%
\ifentrytype{article}{}{\printtext{\bibstring{in}\space}}}
\renewcommand*{\labelalphaothers}{\textsuperscript{+}}
\DefineBibliographyStrings{english}{%
bibliography = {References},
}

\begin{document}

\title{Allocation of Fungible Resources via a Fast, Scalable Price Discovery Method}
\author{Akshay Agrawal \and Stephen Boyd \and Deepak Narayanan \and
Fiodar Kazhamiaka \and Matei Zaharia}

\maketitle

\begin{abstract}
We consider the problem of assigning or allocating resources to
a set of jobs.  We consider the case when the resources are fungible,
that is, the job can be done with any mix of the resources, but
with different efficiencies.  In our formulation we maximize a 
total utility subject to a given limit on the resource usage,
which is a convex optimization problem and so is tractable.
In this paper we develop a custom, parallelizable algorithm for solving the 
resource allocation problem that scales to large problems, with millions of
jobs. Our algorithm is based
on the dual problem, in which the dual variables associated with the resource
usage limit can be interpreted as resource prices. Our method updates the
resource prices in each iteration, ultimately discovering the optimal resource
prices, from which an optimal allocation is obtained. We provide an
open-source implementation of our method, which can solve problems
with millions of jobs in a few seconds on CPU, and under a second
on a GPU; our software can solve smaller problems in milliseconds. On
large problems, our implementation is up to three orders of magnitude
faster than a commerical solver for convex optimization.
\end{abstract}

\section{Introduction}\label{s-res-alloc}

We consider the problem of allocating fungible resources
to a set of jobs. The goal is to maximize a concave utility function of the
allocation, given limits on the amount of available resources. This is a convex
optimization problem, and so is tractable.

For this problem we develop a custom, efficient method amenable to parallel
computation, allowing it to scale to problem sizes larger than can be handled
by off-the-shelf solvers for convex optimization. Our method solves the dual
problem, adjusting the dual variable for the resource constraint to its
optimal value. For a given dual variable value, the dual function splits into
several small resource allocation problems, one per job, which can be solved in
parallel using an analytical solution that we derive. (In this sense our method
can be interpreted as a simple dual decomposition method
\cite[\S6.4]{bertsekas1999nonlinear} \cite[\S3.2]{boyd2007notes}, with
an efficient method for evaluating the dual function.)
Because this dual variable can be interpreted as resource prices
\cite[\S5.4.4]{bv2004convex}, our method has a natural interpretation. Roughly
speaking, each job determines its resource usage independently. Our method
iteratively adjusts the prices to their optimal values, \ie, it discovers the
optimal resource prices. From these, we obtain an optimal allocation.

Our motivating application comes from computer systems. Here, the jobs are
computational tasks that are to be scheduled on a number of interchangeable
hardware configurations (for example, as in \cite{narayanan2020heterogeneity},
where each resource is a different type of GPU). Each allocation or schedule
leads to an estimated throughput, and the quality of the allocation is judged
by a utility function of the achieved throughput. (We note that the resource
allocation problem studied in this paper arises in several other contexts, and
that our method is generically applicable across all of them.)

\paragraph{Outline.} 
We state the resource allocation problem in \S\ref{s-res-alloc}. The remainder
of the paper develops and demonstrates our
price discovery algorithm. In \S\ref{s-duality} we describe the (partial)
Lagrangian, dual function, and dual problem for the resource allocation
problem, and we explain how the dual function can be evaluated, and how an
optimal resource allocation can be found from the optimal dual variables
(prices). In \S\ref{s-subprob} we give an analytical solution to the
subproblems that arise for each job when evaluating the dual function. In
\S\ref{s-price-opt} we give our price discovery algorithm. In
\S\ref{s-examples} we describe our software implementation of the method, which
heavily exploits the parallelism inherent to solving the subproblems, and can
be run on a CPU or a GPU. In this same section we demonstrate our
implementation on some numerical examples, and show that it is often orders of
magnitude faster than a commerical solver for convex optimization. Finally
in \S\ref{s-conclusions} we explain how our problem connects to other types of
resource allocation problems, and mention some extensions to the problem that
are compatible with our method.

\subsection{Related work}

\paragraph{Price discovery methods.} Resource allocation problems
arise in many fields, and they are frequently solved by price adjustment
methods that are similar in spirit to ours. Price adjustment methods are an
instance of a general family of methods called dual decomposition
\cite[\S6.4]{bertsekas1999nonlinear} \cite[\S3.2]{boyd2007notes}, in which
Lagrange multipliers are introduced for complicating constraints in a way that
makes it efficient to evaluate the dual function (and obtain a subgradient).
These methods have been applied widely, especially in communication networks
\cite{kelly1998rate, xiao2004simultaneous, yu2006dual, palomar2006tutorial,
bertsekas1992data} and energy management \cite{finardi2006solving,
zhang2013robust, hu2018distributed}, but also in other contexts 
\cite{yu2006dual, rantzer2009dynamic, komodakis2007mrf, schutz2009supply}. In
communications, it has been shown that under certain conditions, the
TCP/IP protocol can be interpreted as a distributed dual method for solving a
utility maximization problem, with different congestion control
mechanisms optimizing for different utility functions \cite{chiang2007layering}.

\paragraph{Real-time optimization.}
In this paper we develop an extremely fast method for solving a specific
class of convex optimization problems that scales to very large problems (with
tens or hundreds of millions of jobs). Because our method is so fast, it could
conceivably be deployed in a real-time setting, in which the problem would
be solved several times a second, resulting in a new allocation each time
(in the setting of computer systems, this might be reasonable for time-slicing
threads across CPU cores, but less so for moving whole tasks across different
servers). (As we will discuss later, our method also has other uses, such as
pricing resources in a shared or cloud data center.) There is a large body of
work on real-time optimization, for more general classes
of problems than ours. Small to medium-size problems can be solved extremely
quickly using embedded solvers \citep{domahidi2013ecos,osqp,wang2010fast} or
code generation tools that emit solvers specialized to parametric problems
~\citep{Mattingley:2012,chu2013code,osqp_codegen}. For example, the aerospace
and space transportation company SpaceX uses the quadratic program code
generation tool CVXGEN~\cite{Mattingley:2012} to land its
rockets~\citep{spacex}.

Indeed, for slower rates, in which a problem needs to be solved just once
every few minutes, even high-level domain-specific languages for optimization
such as CVXPY \cite{diamond2016cvxpy, agrawal2018rewriting} have been found to
be sufficiently fast, especially when symbolic parameters are used, which
make recompilations of a single problem with different numerical data essentially
free \cite{agrawal2019differentiable}. For example, the technology and media
company Netflix partially replaced the Linux CFS scheduler with a combinatorial
optimization subroutine, implemented using CVXPY, to allocate containers to
CPUs in a way that minimizes interference \cite{rostykus2019predictive}.

\section{Resource allocation problem}\label{s-res-alloc}
In this section we state the resource allocation problem and study some of its
basic properties. In \ref{s-res-alloc-to-jobs}, we lay out the main parts of
the resource allocation problem and introduce the concept of throughput, which
is a linear function of a job's resource allocation. In \ref{s-utility}, we
introduce the concept of utilities, which are functions of the throughput that
measure the quality of an allocation; we also give some examples of utility
functions. In \ref{s-utility-interp}, we give some additional interpretations
of utility functions. Finally in \ref{s-res-alloc-prob}, we tie together these
concepts and present the resource allocation problem in its entirety.

\subsection{Resource allocation to jobs}\label{s-res-alloc-to-jobs}
We consider a setting with $n$ jobs (or processes or tasks), 
labeled $i=1,\ldots, n$, and $m$ types of resources, labeled $j=1,\ldots, m$.
In the problems we are interested in, $n$ is typically large, and 
$m$ is typically small (though the amount of resources available
for each type may be large). We let $x_i \in \reals_+^m$ denote the allocation
of the $m$ resources to job $i$.
We collect these resource allocation vectors into a matrix $X\in \reals^{n \times m}$,
with $i$th row $x_i^T$.
We interpret $X_{ij} = (x_i)_j$ as the fraction of time job $i$ gets to use 
resource $j$.
Thus we have $\ones^T x_i \leq 1$ for each $i$, or matrix terms, $X\ones \leq \ones$,
where $\ones$ is the vector with all entries one and the inequality is elementwise.
We refer to $X$, or the collection of vectors $x_i$, as the resource allocation.

\paragraph{Total resource usage limit.}
The $m$-vector $r = \sum_{i=1}^n x_i = X^T \ones$ gives the total 
usage of each of the $m$ resources.  The total resource usage
cannot exceed a given limit $R \in \reals_+^m$, \ie, $r \leq R$.
(We mention that we can easily handle the case in which some jobs consume
more than one unit of resource while running, in which case the constraint
becomes $X^T d \leq R$, where $d \in \reals^m_{++}$ gives the amount of resources
demanded by each job.)

\paragraph{Throughput.}
The throughput of job $i$ is $t_i = a_i^Tx_i$, where $a_i \in \reals_+^m$ is
a given efficiency vector.
The particular form $t_i=a_i^Tx_i$ says that job $i$ can be carried out using any
mixture of the resources, with
$(a_i)_j$ interpreted as the effectiveness or efficiency of using resource $j$ for 
job $i$.
Another interpretation is that the resources are fungible, \ie, they can be 
substituted for each other. We obtain the same throughput $t_i$ for any 
allocation that satisfies $a_i^Tx_i =t_i$.  In particular,
we can `exchange' resource $j$ for resource 
$j'$, by decreasing $(x_i)_j$ by $\delta>0$, and increasing $(x_i)_{j'}$
by $(a_i)_j/(a_i)_{j'} \delta$
(assuming these changes do not violate the constraints $x_i \geq 0$,
$\ones^T x_i \leq 1$).
We can interpret 
$(a_i)_j/(a_i)_{j'}$ as the exchange rate between resource $j$ and $j'$,
for job $i$.

Note that the throughput of job $i$ ranges between the minimum value $t_i=0$
(obtained with $x_i=0$) and a maximum value $t_i = \max_j (a_i)_j$, 
obtained with $x_i = e_q$, where $q = \argmax_j (a_i)_j$.
(In other words, the maximum throughput for a given job 
is obtained by using the most efficient resource, at $100\%$.)

\subsection{Utility}\label{s-utility}
The utility of the allocation to job $i$ is given by $u_i(t_i)$,
where $u_i: \reals_{++} \to \reals$ is the utility associated 
with job $i$, for $i=1, \ldots, n$.  
We will assume these are nondecreasing and concave functions.
Nondecreasing means that we derive more (or the same) 
utility from higher throughput,
and concavity means that there is 
decreasing marginal utility as we increase the throughput.
The total utility is given by $U(t) = \sum_{i=1}^n u_i(t_i)$,
where $t \in \reals_+^n$ is the vector of job throughputs.
The average utility, which can be more interpretable than the total
utility, is $U(t)/n$.

Below, we give a few examples of utility functions.

\paragraph{Linear utility.}
The simplest utility is linear utility, with $u_i(t_i)=t_i$;
in this case the overall utility is the total throughput,
and the average utility is the average throughput.
Roughly speaking the linear utility gives equal weight to 
increasing throughput; nonlinear concave utilities give
more weight to increasing the throughput of a job when the
throughput is small.

\paragraph{Worst-case or min utility.}
A utility function that is used in some applications is the minimum 
throughput or worst-case utility $U(t) = \min_i t_i$.
This utility function is not separable, and so does not fit our
requirement of separability.  Nevertheless we will see below that it is 
can be approximated by separable utilities.
We note that the min-utility is at the opposite extreme from
the linear utility, since roughly speaking it gives no weight to 
increasing any utility above the minimum, and focuses all its
attention on the jobs with minimum throughput.

\paragraph{Log utilities.}
A commonly used strictly concave utility function is the
logarithmic utility 
\BEQ\label{e-log-util}
u_i(t_i) = \log t_i,
\EEQ
which is used in economics (\eg, in Kelly gambling \cite{kelly1956new,
maclean2011kelly, busseti2016risk}) and networking, where it leads to
allocations that are called \emph{proportionally fair} \cite{kelly1998rate}.

\paragraph{Power utilities.}
Another family of strictly concave utility functions is
the power utility $u_i(t_i) = t_i^p$, with $p \in (0,1]$ 
or $u_i(t_i) = -t_i^p$ with $p<0$.
These utility functions are widely used in economics,
where they are called the constant
relative risk aversion (CRRA) or isoelastic utilities
\cite[\S1.7]{eeckhoudt2011economic}. For $p$ positive and small, or negative
and large, the power utility approximates the min-throughput utility (up to a
constant), since it gives much higher weights to smaller throughputs than
larger throughputs.

Log and power utility functions are sometimes described as one family of
utility functions, called
$\alpha$-fairness \cite{mo2000fair},
with the form
\[
  u_i(t_i) = \begin{cases}
    \frac{1}{1 - \alpha}t_i^{1 - \alpha} & \alpha \geq 0 \text{ and } \alpha \neq 1 \\
    \log t_i & \alpha = 1.
  \end{cases}
\]
The choice $\alpha = 0$ yields linear utility, while $\alpha = 1$ yields the log
utility. In networking, it has been shown that taking $\alpha \to \infty$
yields \emph{max-min fairness} \cite{mo2000fair} (in practice, a large value of
$\alpha$ suffices).

\paragraph{Target-priority utility.}
Another useful family of utility functions is based on a target
throughput and a priority,
\BEQ\label{e-target-priority-util}
u_i(t_i) = w_i \min\{ t_i-t_i^\text{des} , 0\},
\EEQ
where $w_i$ is a positive weight parameter
and $t_i^\text{des}$ is a positive target throughput.
This utility is zero when the throughput meets or exceeds the target value,
and decreases linearly, with slope $w_i$, when the throughput
comes short of the target.
The parameter $w_i$ encodes the priority of job $i$, with
higher weight giving higher priority.
With target-priority utility, the total utility is zero if all
job target throughputs are achieved, and negative otherwise;
it is the total of (weighted) shortfalls.
These utility functions are not differentiable, or strictly increasing,
or strictly concave.

\subsection{Utility interpretations}\label{s-utility-interp}
\paragraph{Utility-derived averages.}
Utility functions, and the resulting utility, are meant to measure
the quality of an aggregate throughput.  Linear utility treats all throughputs,
large and small, the same; concavity or curvature of a utility function puts
more weight on the smaller job throughputs than larger ones.
(The extreme here is the worst-case utility, which focuses all its attention
on the smallest throughput.)
When the same utility $u$ is used for all jobs, and $u$ is invertible,
we can interpret the quantity
\BEQ\label{e-u-avg}
u^{-1}\left(\frac{1}{n}\sum_{i=1}^n u(t_i)\right)
\EEQ
as a kind of average of the throughputs, skewed toward the smaller ones.
It has the same units and scale as the throughput itself, and coincides with
well known averages for some choices of utilities.  For example,
it is the (arithmetic) average for linear utility, the geometric mean
for log utility, and the harmonic mean for the inverse utility $u(t_i)=-1/t_i$.
The latter two have been proposed as measures of system performance
that in some cases are more appropriate than the simple arithmetic average
\cite[\S1.8]{hennessy2011computer}.

\paragraph{Connection to risk-adjusted average throughput.}
Utilities are closely related to the concept of risk-adjusted average
throughput.
Let $\avg(t)$ denote the average throughput and $\var(t)$ denote the variance of 
the throughput across jobs, \ie,
\[
\avg(t) = \frac{1}{n}\sum_{i=1}^n t_i, \qquad
\var(t) = \frac{1}{n}\sum_{i=1}^n t_i^2 - \left(
\frac{1}{n}\sum_{i=1}^n t_i \right)^2.
\]
The average throughput is a natural measure of overall throughput;
the variance is a natural measure of fairness since it quantifies
how different the job throughputs are.
The risk-adjusted throughput is defined as
\[
\avg(t) - \frac{\gamma}{2} \var(t), 
\]
where $\gamma>0$ is the so-called risk aversion parameter.
It measures an aggregate throughput, with an adjustment for fairness,
scaled by $\gamma$.
The risk-adjusted throughput metric is large when the 
average throughput is large and 
the variation in throughput across the jobs is small.
If we maximize it, it means
we will accept a reduction in the average throughput,
if it comes with a sufficient decrease in the variance of the throughputs.
This concept is widely used in finance, especially in portfolio construction,
where it dates back to the 1950s \cite{markowitz1952portfolio, tobin1965theory}.

Now suppose that $\phi:\reals_+ \to \reals$ is concave, increasing, and 
twice differentiable, with
$\phi(0)=0$, $\phi'(0)=1$, and $\phi''(0)=-1$.
(For example, $\phi(a) = 1-e^{-a}$.)
We can define a family of utility functions $u_i(t_i) = \phi(\gamma t_i)$, 
where $\gamma>0$.  A basic result is that small $\gamma$, we have 
\[
\frac{1}{\gamma} \phi^{-1} (U(t)/n) =
\avg(t) - \frac{\gamma}{2} \var(t) + o(\gamma^2).
\]
(This can be seen by taking second-order Taylor expansions of
the utility functions $u_i$ and $\phi^{-1}$ around $0$, and dropping
higher order terms.)
The lefthand side is the utility average \eqref{e-u-avg} defined by $u_i$,
scaled by $1/\gamma$; the righthand side is the risk-adjusted average
throughput, plus higher order terms in $\gamma$.
Thus for small $\gamma$, the utility, mapped through the monotone
increasing mapping $a \mapsto (1/\gamma) \phi^{-1}(a/n)$,
is approximately the risk-adjusted average throughput.

Connections between utility maximization and risk aversion have been
studied extensively in economics; \eg, 
see \cite{friedman1948utility, tobin1958liquidity, pratt1964risk,
arrow1971essays, gollier2001economics}.

\subsection{Resource allocation problem}\label{s-res-alloc-prob}
The resource allocation problem is
\BEQ\label{e-res-alloc}
\begin{array}{ll}
\mbox{maximize} & U(t)\\
\mbox{subject to} & x_i \geq 0, \quad \ones^T x_i \leq 1, \quad t_i=a_i^Tx_i,
\quad i=1, \ldots, n\\
& \sum_{i=1}^n x_i \leq R,
\end{array}
\EEQ
with variables $x_i \in \reals^m$ and $t_i \in \reals$, $i=1, \ldots, n$.
The problem data are the utility functions $u_1, \ldots, u_n$, the 
efficiency vectors $a_1, \ldots, a_n$, and the total resource usage limit $R$.
We recall that $x_i$ is the allocation of resources to job $i$, $t_i = a_i^T x_i$ is
the throughput achieved by job $i$, the constraint $\ones^T x_i \leq 1$ says
that each job consumes resources for at most the total schedulable time, and
the constraint $\sum_{i=1}^{n} x_i \leq R$ says that jobs cannot consume more
resources than there are at hand. We denote an optimal allocation as
$x_i^\star$, $i=1,\ldots, n$, and its associated optimal throughput as
$t^\star$, and utility as $U^\star = U(t^\star)$.

The resource allocation problem \eqref{e-res-alloc} is a convex optimization problem
and therefore tractable \cite[\S1]{bv2004convex}.
We observe that the objective and the first line of constraints
are separable across jobs; the total resource usage limit (the last constraint)
couples the different jobs.
In other words, without the last total resource usage constraint, the resource
allocation problem splits into $n$ separate problems, one for each job $i$.
Our method will leverage this idea.

The resource allocation problem is infeasible only in obvious pathological
cases such as $a_i=0$ with log utility.  We assume henceforth that the problem is 
feasible.  It always has a solution, since the feasible set is compact.
The problem is bounded above; the maximum possible utility is 
$\sum_i u_i(\max_j (a_i)_j)$, the utility when all throughputs take on
their maximum possible values.
Our analysis later will show that the solution need not be unique, 
even when the utilities are strictly concave.

In the resource allocation problem \eqref{e-res-alloc} we can replace
the utility $U(t)$ with the average utility $U(t)/n$ or, when the utilities
$u_i$ are the same and invertible, the utility average \eqref{e-u-avg}.
These monotone transformations of the utility yield an equivalent problem.

\paragraph{Pareto optimality.} If the utility functions are strictly
increasing, then the throughput $t$ achieved by an optimal allocation $X$
is Pareto efficient \cite[\S4.7.3]{bv2004convex}, in the following sense.
If $\tilde t$ is a throughput vector for another feasible allocation for which
$\tilde t_i > t_i$ for some job $i$, then there must be some other job $j$ for
which $\tilde t_j < t_j$ (if this were not the case, then evidently $X$
would not be optimal).

\paragraph{Application to computer systems.} 
In our motivating application,
the jobs represent computational tasks or services.
The resources represent
different hardware configurations on which the jobs may be run, 
such as types of CPUs (differing in cache sizes, clock frequency, 
core count), GPUs (differing in memory, core count); 
or servers (differing in, say, CPU, GPU, RAM, storage,
network bandwidth). The resources are fungible because a job can be run on any
of the $m$ hardware configurations, but with different efficiencies. The entry
$(x_i)_j$ is the fraction of time job $i$ will run on hardware type $j$, and
$(a_i)_j$ is the throughput job $i$ would obtain if it were to run entirely on
hardware $j$ (the entries of $a_i$ can be obtained by profiling each job on the
different hardware configurations). Throughput can be measured in several ways,
such as number of individual tasks that can be processed per second, or the
number of MIPS (millions of instructions per second) that a hardware
configuration can achieve on a job type.

We can interpret the resource allocation problem as finding an optimal
time-slicing of the $n$ jobs across the $m$ hardware configurations, where
optimality is measured by total utility. This problem (with linear and
worst-case utility) was recently studied in \cite{narayanan2020heterogeneity},
as part of a system for scheduling deep learning jobs on GPUs. (The problem of
scheduling tasks across interchangeable servers was studied in
\cite{tumanov2016tetrisched}, but using a different formulation of the
allocation problem than ours.)

In addition to finding an optimal time-slicing of the jobs, our method
discovers the optimal prices of the different hardware configurations.
These prices could be used to inform markets that provide
several users access to a shared pool of computational resources, as they
convey the value of each resource relative to the others. They could also be
used to actually charge jobs for hardware usage.

\section{Duality and resource prices}\label{s-duality}
In the remainder
of the paper, we
develop and demonstrate our price discovery method for efficiently
solving \eqref{e-res-alloc}. We recall that our method is based on solving
the dual problem: we introduce a Lagrange multiplier for the constraint
on resource usage, which splits the dual function into one small
resource allocation problem per job. These subproblems can be solved
efficiently and in parallel using an analytical solution. This lets us evaluate
the dual function and a subgradient very cheaply; we use the subgradients to
adjust the prices to their optimal values, from which we obtain an optimal
allocation.

In this section we cover some standard results about duality in convex 
analysis and optimization that we will use in our solution method.
For background on duality in convex optimization,
see \cite[Chap.~5]{bv2004convex}, \cite[Chap.~XII]{urruty1993convex}, or
\cite[\S VI.28]{rockafellar1970convex}.

\paragraph{Lagrangian and dual function.}
We first reformulate the problem \eqref{e-res-alloc} as 
\[
\begin{array}{ll} \mbox{maximize} & U(t) - \mathcal I(X,t) \\
\mbox{subject to} &r  \leq R,
\end{array}
\]
with variables $X$ and $t$, where $r=X^T\ones$ is the total resource
usage, and $\mathcal I$ is the indicator function of the constraints
\[
x_i \geq 0, \quad  \ones^T x_i \leq 1, \quad t_i = a_i^T x_i, \quad i=1, \ldots, n.
\]
(This means that $\mathcal I(X,t)=0$ when these constraints are
satisfied, and $\mathcal I(X,t)=\infty$ when they are not.)
We introduce $p\in \reals_+^m$ as a dual variable 
(or Lagrange multiplier or shadow price)
for the constraint $r \leq R$, and form the Lagrangian
\[
L(X,t,p) = U(t) - \mathcal I(X,t) -p^T (r-R).
\]
We can interpret $p$ as a set of prices for the resources
\cite[\S5.4.4]{bv2004convex}, and the part of the Lagrangian $U(t)-p^Tr$ as the
net utility, 
\ie, the utility derived from the throughput minus the cost of using
the resources, at the prices given by $p$.

The dual function is defined as
\[
g(p) = \max_{X,t} L(X,t,p).
\]
The dual function is convex.
This is the optimal value of the Lagrangian for the resource price vector $p$.

\paragraph{Evaluating the dual function.}
We first observe that the Lagrangian is separable across jobs $i$, \ie,
a sum of functions of $x_i$ and $t_i$:
\BEAS
L(X,t,p) &=& 
U(t)-p^T (r-R)  - \mathcal I(X,t) \\
&=& p^T R + \sum_{i=1}^n \left(u_i(t_i) - p^Tx_i - \mathcal I(x_i \geq 0,~
\ones^T x_i \leq 1,~t_i=a_i^Tx_i)   \right).
\EEAS
To evaluate the dual function $g(p)$ we maximize this over all $x_i$ and $t_i$;
by separability we can maximize separately for each $i$.
For $i=1, \ldots, n$ we solve the problem
\BEQ\label{e-subprob}
\begin{array}{ll} \mbox{maximize} & u_i(t_i) - p^Tx_i \\
\mbox{subject to} & x_i \geq 0, \quad \ones^Tx_i \leq 1, \quad t_i=a_i^Tx_i.
\end{array}
\EEQ
This is a small convex optimization problem with $m+1$ variables,
which we will show to how to solve analytically in \S\ref{s-subprob}.
Its objective is the net utility for job $i$, \ie, the utility minus 
the cost of resources used.
The dual function value $g(p)$ is the sum of the optimal values of the
subproblems \eqref{e-subprob}, plus $p^TR$.

\paragraph{Optimal value bounds from the dual function.}
Since for any feasible $X,t$ and any $p \in \reals_+^m$
we have $L(X,t,p) \geq U(t)$, it follows that
\BEQ\label{e-upper-bnd}
g(p) \geq U^\star.
\EEQ
In other words, the dual function gives an upper bound on the optimal
utility of the resource allocation problem \eqref{e-res-alloc}.

\paragraph{Dual problem.}
The dual problem has the form
\BEQ\label{e-dual}
\begin{array}{ll} \mbox{minimize} & g(p) \\
\mbox{subject to} & p \geq 0,
\end{array}
\EEQ
with variable $p \in \reals^m$.
This has the natural interpretation of choosing
$p$ to obtain the best (\ie, smallest) 
upper bound on $U^\star$ in \eqref{e-upper-bnd}.
The dual problem \eqref{e-dual} is convex.  
We denote an optimal $p$ as $p^\star$, and refer to $p^\star$ as 
the \emph{optimal resource prices}.
Solving the dual problem is sometimes called \emph{price discovery}, since 
solving it finds the optimal prices.

\paragraph{Strong duality.}
A standard result from convex optimization states that
\[
g(p^\star) = U^\star,
\]
\ie, the upper bound on optimal utility from the dual function \eqref{e-upper-bnd}
is tight when evaluated at the optimal prices. 
(Strong duality holds here since the only constraints are linear equalities and 
inequalities; see \cite[\S 5.2.3]{bv2004convex}.)

\paragraph{Recovering an optimal allocation from optimal prices.}
A basic duality result from convex optimization states that any optimal
allocation $X^\star$ and throughput $t^\star$ maximizes
$L(X,t,p^\star)$ over $X$ and $t$. But since $L$ is separable across jobs $i$,
this means that for each job $i$, $x_i^\star$ and $t_i^\star$ are solutions of 
the problem
\[
\begin{array}{ll} \mbox{maximize} & u_i(t_i) - (p^\star)^Tx_i \\
\mbox{subject to} & x_i \geq 0, \quad \ones^Tx_i \leq 1, \quad t_i=a_i^Tx_i.
\end{array}
\]
This has a very nice interpretation.   Each job derives utility $u_i(t_i)$,
and pays for the resources consumed at the optimal prices $p^\star$.  The optimal
allocation maximizes the net utility, \ie, the utility derived from the resources
minus the amount paid for the resources.
This relation between optimal allocation and optimal prices is
key to our price-based method.

\paragraph{Resource prices.}
The interpretation of $p^\star$ as a set of resource prices is 
standard throughout applications that use optimization.
To explain the interpretation, define $U(R)$ to be the optimal value of
the resource allocation problem \eqref{e-res-alloc}, \ie, the maximum 
utility, as a function of the total resource usage limit $R$.
A standard duality result is
\BEQ\label{e-pstar-grad}
p^\star = \nabla_R U^\star(R),
\EEQ
the gradient of the optimal utility with respect to the resource limits.
Thus $p_j^\star$ is the (approximate) increase in optimal utility obtained per 
increase in resource $j$.
(When $U^\star(R)$ is not differentiable, we replace
the gradient above with a subgradient of $-U$.)
The interpretation of the partial derivative of maximum utility with respect
to a resource limit as a price for the resource is common in many
fields, \eg, in communication networks \cite{kelly1956new} and power networks
\cite{bohn1984optimal}.

\paragraph{Subgradient of dual function.}
We mention for future use how to find a subgradient of $g$ at $p$.
Let $\tilde x_i$ and $\tilde t_i$ be the solutions of the subproblems
\eqref{e-subprob}, for $i=1, \ldots, n$.
Then a subgradient of $g$ is given by
\BEQ\label{e-g-subgrad}
q = R-r = \sum_{i=1}^n R - \tilde x_i
\EEQ
If $g$ is differentiable at $p$, then $q = \nabla g(p)$.
The subgradient $q$ has a nice interpretation: it is the difference between
the total resource limit $R$ and
the total resource usage $r$, when you choose the allocations via the subproblems
\eqref{e-subprob} with resource prices $p$.

\section{Solving the subproblem} \label{s-subprob}
In this section we explain how to analytically solve the subproblems \eqref{e-subprob}.
In this section we will drop the subscript $i$ (which indexes the jobs), to keep
the notation light, and express  the problem as
\BEQ\label{e-subproblem}
\begin{array}{ll} \mbox{maximize} & u(t) - p^Tx \\
\mbox{subject to} & x \geq 0, \quad \ones^Tx \leq 1, \quad t=a^Tx,
\end{array}
\EEQ
with variables $x\in \reals^m$ and $t \in \reals$.
We write this as
\BEQ\label{e-opt-over-t}
\begin{array}{ll} \mbox{maximize} & u(t) - c(t), \end{array}
\EEQ
with variable $t \in \reals$, 
where $c(t)$ is the optimal value of the linear program (LP)
\[
\begin{array}{ll} \mbox{minimize} & p^Tx \\
\mbox{subject to} & x \geq 0, \quad \ones^Tx \leq 1, \quad a^Tx = t,
\end{array}
\]
with variables $x\in \reals^{m}$ and $t\in \reals$.
We now show how to solve this LP analytically.  Indeed, we will give the solution
parametrically, and obtain an explicit formula for $c(t)$.

\subsection{Parametric solution of the LP}
We introduce a slack variable $s\in \reals$ and express it in standard form
\[
\begin{array}{ll} \mbox{minimize} & p^Tx \\
\mbox{subject to} & (x,s) \geq 0\\
& \left[ \begin{array}{cc} a^T & 0 \\ \ones^T & 1 \end{array} \right] 
\left[ \begin{array}{c} x \\  s \end{array} \right] =
\left[ \begin{array}{c} t \\  1 \end{array} \right].
\end{array}
\]
A basic result for LPs states that there is always a basic feasible solution,
\ie, one in which at most two entries of $(x,s)$ are nonzero
\cite[\S 2.2]{bertsimas1997introduction}. (Two is the 
number of linear equality constraints.) There are $m(m+1)/2$ such choices of
two nonzero entries of $(x,s)$. This tells us that the resource allocation
subproblem always has a solution that uses at most two of the resources. 

\begin{figure}
\includegraphics{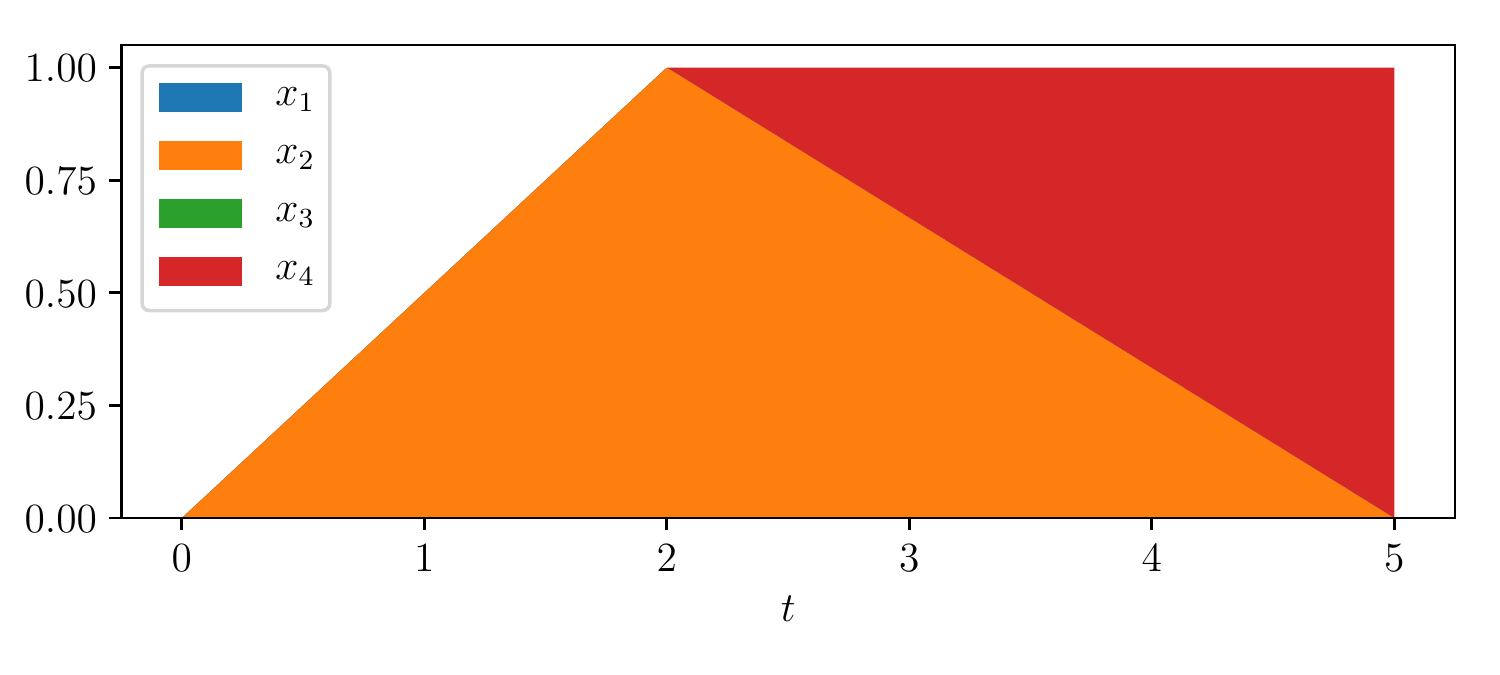}
\caption{Basic feasible solutions for a subproblem with $m=4$ resources,
with data \eqref{e-illus},
as $t$ varies.  The first and third resources are never used.}
\label{f-basic-feasible-solution}
\end{figure}

This is illustrated in figure~\ref{f-basic-feasible-solution}
for a subproblem with $m=4$, and
\BEQ\label{e-illus}
a = (1,2,3,5), \qquad
p = (1,1,4,6).
\EEQ
The figure shows a basic feasible solution of the subproblem as $t$
ranges from $0$ to $5$, its range of feasible values.
For $0 \leq t \leq 2$, only the second 
resource (which has the highest value of $a_j/p_j$) is used;
for $2 < t < 5$, both the second and fourth resources are used;
and for $t=5$, only the fourth resource 
(which has the largest value $a_j$) is used.
We note that the first and third resources are never used.

We will now assume that the efficiencies are sorted and distinct, so 
$a_1 < \cdots < a_m$.
(The results given below are readily extended to the case when they are not 
distinct.)
First suppose that $x_i$ and $x_j$ are  nonzero, with $i <j$.
(This corresponds to the case when the job uses only resources $i$ and $j$.)
Then $s=0$, so $x_i+x_j=1$ and $a_ix_i + a_jx_j = t$, so
\BEQ\label{e-xixj}
x_i = \frac{a_j-t}{a_j-a_i}, \qquad
x_j = \frac{t-a_i}{a_j-a_i}.
\EEQ
For these to be nonnegative, we must have $t\in[a_i,a_j]$.
The associated objective value is 
\BEQ\label{e-optcostij}
p_i \frac{a_j-t}{a_j-a_i} + p_j \frac{t-a_i}{a_j-a_i}.
\EEQ
As $t$ varies between $a_i$ and $a_j$, this varies affinely
between $p_i$ and $p_j$, respectively.

Now consider the special case when $x_j$ and $s$ are nonzero.
(This corresponds to the case when the job uses only resource $j$.)
Then we have $a_jx_j = t$, so
\[
x_j = t/a_j.
\]
For $x_j$ and $s=1-x_j$ to be nonnegative, we need $t \in [0,a_j]$.
The corresponding objective value is
\[
p_j t/a_j.
\]
These are the same as the formulas \eqref{e-xixj} and \eqref{e-optcostij} above, 
with $i=0$, where we define $a_0=p_0=0$.

The optimal value $c(t)$ is the pointwise 
minimum of the $m(m+1)/2$ affine functions, restricted to an interval,
\[
p_i \frac{a_j-t}{a_j-a_i} + p_j \frac{t-a_i}{a_j-a_i} + \mathcal I(t \in [a_i,a_j]),
\quad i<j.
\]
The graphs of these functions are line segments that
connect the points $(a_i,p_i)$ and $(a_j,p_j)$, including the additional point 
$(a_0,p_0)= (0,0)$.
The graph of $c(t)$ is the pointwise minimum of these.  
It is easy to see that $c$ is piecewise affine with kinks points that 
are a subset of the values $a_1, \ldots, a_m$.
It is increasing and convex, and satisfies $c(0)=0$, and has domain $[0,a_m]$.

\begin{figure}
\includegraphics{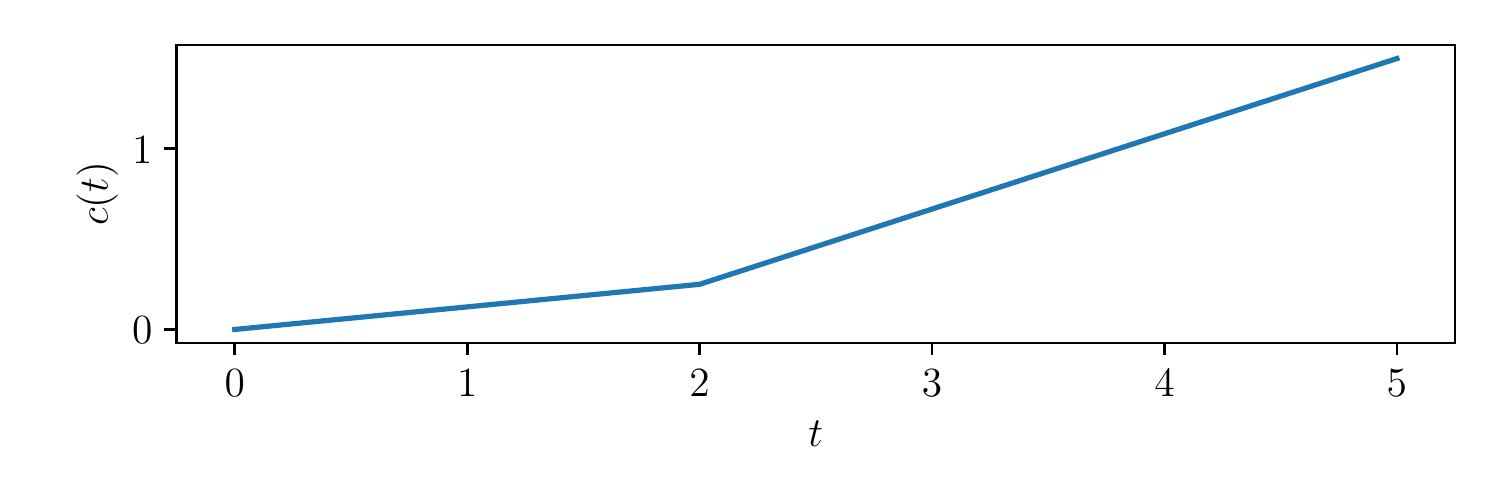}
\caption{A graph of the piecewise affine optimal cost $c(t)$, for a fixed price
vector and four resources.}
\label{f-cost-versus-throughput}
\end{figure}

Given $a$ and $p$, it is straightforward to directly compute the piecewise 
affine function $c$, \ie, to find its kink points and the slope in between
successive kink points.
It is completely specified by a subset of $\{0,\ldots,m\}$, 
given by $0=i_1< \cdots  < i_r = m$.
The kink points are $a_{i_l}$, $l=1,\ldots, r$, and the value of $c$ at these
points is $p_{i_l}$. 

This is illustrated in figure~\ref{f-cost-versus-throughput},
for the same example as above with data \eqref{e-illus}.
In this case the active subset is
$0,2,4$, and $c$ is piecewise affine with kink points $0,a_2=2,a_4=5$,
and associated values $0,p_2=1,p_4=6$.
\begin{figure}
\includegraphics{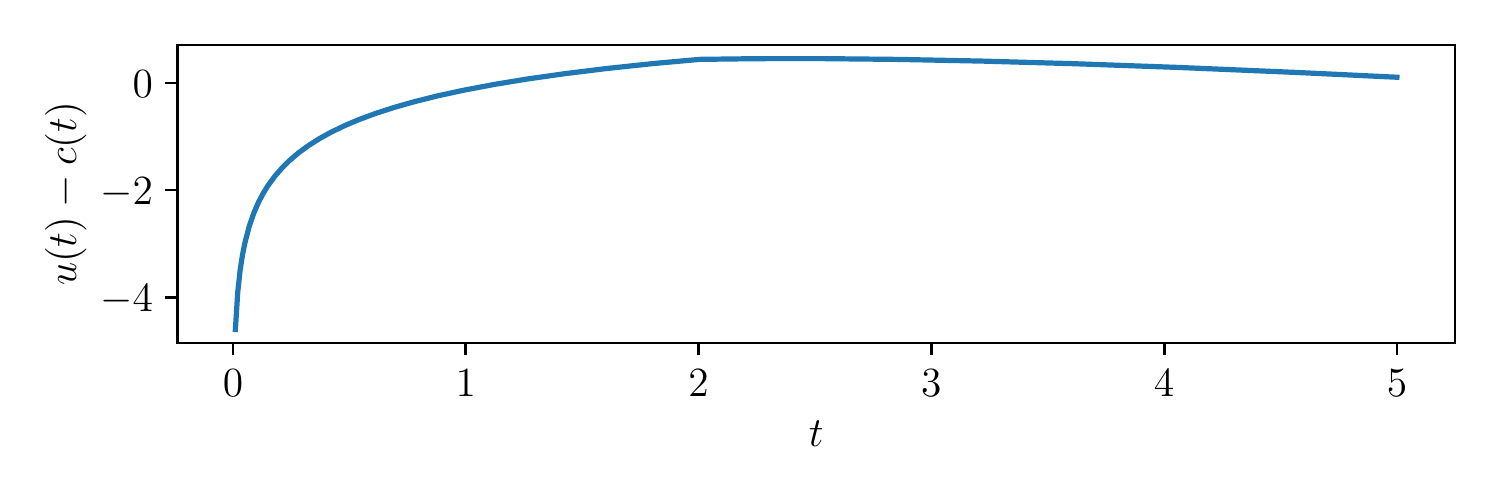}
\caption{Net utility, \ie, utility minus cost,
for a log utility and four resources.}
\label{f-excess-utility}
\end{figure}

\subsection{Solving the subproblem}
Once we have the explicit function $c$, we can readily solve the 
one dimensional problem \eqref{e-opt-over-t}.
To do this we maximize $u$ minus an affine function over each of the intervals 
between successive kink points and choose the one with largest objective.
Thus we need to maximize 
$u(t) - (\alpha t + \beta)$ over $t \in [\gamma,\delta]$, where $\alpha,
\beta, \gamma, \delta$ are given, with $\alpha >0$.
This is readily done for any utility function.
For example, if $u$ is increasing and strictly concave 
(\eg, log utility) we have the explicit formula
\BEQ\label{e-opt-t-log}
t = \Pi (u')^{-1}(\alpha),
\EEQ
where $\Pi$ is the projection onto the interval $[\gamma,\delta]$.
(Since $u$ is strictly concave, $u'$ is increasing and therefore invertible.)
Figure~\ref{f-excess-utility} plots $u(t) - c(t)$ for our running example
with $m=4$ resources, data \eqref{e-illus}, and log utility;
in this case, a value of $t \approx 2.4$ is optimal.

This small problem is also easily solved for cases when the utility is
neither increasing nor strictly concave, \eg, the target-priority utility
\eqref{e-target-priority-util}.
For this utility function, the solution is
\BEQ\label{e-opt-t-target}
t = \left\{ \begin{array}{ll} 
t^\text{des} & w>\alpha, ~ \gamma \leq t^\text{des} \leq \delta, \\
\gamma & w \leq \alpha,\\
\delta & w>\alpha, ~ \delta \leq t^\text{des}.
\end{array}\right.
\EEQ

\paragraph{Summary.} We summarize the algorithm for solving the 
job subproblem \eqref{e-target-priority-util} below.

\begin{algdesc}{\sc Maximizing net utility} \label{alg-max-net-utility}

\textbf{given} prices $p \in \reals^m_{+}$, efficiency vector $a \in
\reals^m_{+}$ ($a_1 < a_2 < \cdots < a_m$), utility function $u : \reals \to
\reals$

\begin{enumerate}
\item \emph{Compute $c(t)$.} Compute the piecewise affine function $c$ as
the pointwise minimum of $m(m+1)/2$ affine functions.
\item \emph{Compute optimal throughput.} Compute the $t$ 
maximizing $u(t) - c(t)$ (\eg, via \eqref{e-opt-t-log} or 
\eqref{e-opt-t-target})
\item \emph{Compute optimal allocation.} Compute the allocation for the
pair $(i, j)$ achieving the maximum net utility, as in \eqref{e-xixj}.
\end{enumerate}
\end{algdesc}
Note that we access the utility function $u$ in just two ways:
we must be able to evaluate $u$ at any throughput $t$, and obtain the
$t$ maximizing $u(t) - c(t)$.

\paragraph{Complexity.}
The complexity of this method of solving the job subproblem is quadratic in
$m$. By solving it in parallel for all $n$ jobs (in \S\ref{s-examples}, we
describe one way to parallelize these solves) , we go from a proposed resource
price vector $p$ to an allocation that satisfies all constraints, except
possibly the total resource usage limit. At the same time we can evaluate the
dual function $g(p)$, and a subgradient of it, as described above.

\section{Price discovery algorithm}\label{s-price-opt}

Now we can give our method for solving the resource allocation 
problem \eqref{s-res-alloc}.
We will instead solve the dual problem \eqref{e-dual},
in order to discover the optimal resource prices; we recover
an optimal allocation from solutions of the subproblems \eqref{e-subprob}.

Our dual or price discovery algorithm
adjusts the resource prices.  We let $p^k\in \reals_+^m$ denote the
prices in iteration $k$.
We first evaluate the dual function $g(p^k)$ and a subgradient $q^k$.
We do this by solving the subproblems as in \S\ref{s-subprob} for each
$i$, using algorithm~\ref{alg-max-net-utility}, in parallel.  
This gives us an upper bound on $U^\star$, and a 
resource allocation $X^k$ that satisfies the job constraints
$x_i \geq 0$ and $\ones^T x_i \leq 1$,
but need not satisfy the total resource usage
limit $(X^k)^T \ones \leq R$.

From this allocation we create a feasible allocation, by scaling 
down each column of $X^k$ so that the resource usage limit holds.  
We refer to this feasible allocation as $\tilde X^k$.  We evaluate its utility,
which is a lower bound on $U^\star$.
Thus we have lower and upper bounds on $U^\star$,
\[
U(\tilde X^k) \leq U^\star \leq g(p^k).
\]
We refer to $g(p^k)-U(\tilde X^k)$ as the duality gap in iteration $k$.
We quit when this is small, with a guarantee on how far from optimal the 
current allocation is.

\paragraph{Standard subgradient price update.}
To move to the next iteration we update the prices.
A standard subgradient method uses the projected gradient update
\[
p^{k+1} = \max\{ p^k - \alpha_k q^k, 0\},
\]
where $\alpha_k$ are positive step lengths that satisfy $\alpha_k \to 0$ and 
$\sum_{k=1}^\infty \alpha_k = \infty$.
This simple update guarantees convergence, \ie, $p^k \to p^\star$
\cite{shor1985nondifferentiable, boyd2003subgradient}. It is also extremely
intuitive. Recall that $q^k = R-r^k$, so that $-q^k$ tells us how much we are
over-using the resources, \ie, $-q^k_j > 0$ means that using the prices $p^k$,
we have $r^k_j > R_j$. The subgradient update above says that we should 
increase the price for resources we are currently over-using, and decrease
the price for any resource we are under-using 
(but never decrease a price below zero).

Recall that the optimal price vector gives the true prices of each resource,
from which an optimal allocation is easily obtained (see \S\ref{s-duality}).

\paragraph{More efficient price updates.}
The projected subgradient method described above always works,
even when $g$ is nondifferentiable.  
We can also use more efficient and sophisticated methods
for minimizing $g$ that rely on $g$ being differentiable,
for example a quasi-Newton method such as the BFGS method
\cite{broyden1970convergence, fletcher1970new, goldfarb1970family, shanno1970conditioning}
or its limited memory variants \cite{nocedal1980updating, liu1989limited}.
While $g$ need not be differentiable, for example with target-priority 
utility, it is quite smooth when $n$ is large, and we have observed no
practical cases where it failed.  (In any case, we can always fall back on 
the basic subgradient method.)

\paragraph{Initialization.}
The initial price vector $p^1$ can be anything, including $0$ or $\ones$.
We have found a simple initialization that depends on the 
data ($a_i$, $R$, and $u_i$) and yields prices that are reasonably close
to the optimal ones.
We start with the simple allocation
$x_i = (1/n)R$ for all $i$, which distributes the full usage budget uniformly
across the jobs. (If $\kappa = \ones^T (1/n) R >1$, we take
$x_i = (1/(\kappa n)R$, so $\ones^T x_i = 1$.)
We take as a starting price
\BEQ\label{e-initialization}
p^1 = \nabla_R U(R) = \frac{1}{n}\sum_{i=1}^n u_i'(a_i^Tx_i) a_i.
\EEQ
(When $u_i$ is not differentiable, we can take a supergradient, \ie,
the negative of a subgradient of $-u$.)
This simple initialization is motivated by the observation that the optimal
prices give the sensitivity of the total utility to the resource constraint
\cite[\S5.6.3]{bv2004convex}. We have found it to work well in practice.

\paragraph{Algorithm summary.} We summarize our price discovery algorithm
for solving the resource allocation problem \eqref{e-res-alloc} below.

\begin{algdesc}{\sc Resource allocation price discovery} \label{alg-price-opt}

\textbf{given} efficiency vectors $a_i \in \reals^m_{+}$,
utility functions $u_i : \reals \to \reals$, 
total resource usage limit $R$, initial resource prices $p^1\in \reals_+^m$,
tolerance $\epsilon >0$

For iteration $k=1,2, \ldots, K^\text{max}$

\begin{enumerate}
\item With prices $p^k$, solve $n$ subproblems in parallel
to find 
\begin{itemize}
\item allocation $X^k$
\item dual function value $g(X^k)$
\item dual function gradient $q^k$
\item feasible allocation utility value $U(\tilde X^k)$
\end{itemize}
\item Quit if $g(p^k) - U(\tilde X^k) \leq \epsilon$
\item Update prices to obtain $p^{k+1}$
\end{enumerate}
\end{algdesc}

A reasonable choice of the tolerance $\epsilon$ is $10^{-3}n$, which
guarantees that our final allocation is no more than $10^{-3}$
suboptimal in average utility.
With log utility, this means the throughputs are optimal to within around 
0.1\%, which is far more than good enough for any practical application.
(The algorithm can be run to much higher accuracy as well.)

\section{Numerical examples}\label{s-examples}

\subsection{Implementation}
We have implemented the price discovery algorithm~\ref{alg-price-opt}
in PyTorch \cite{paszke2019pytorch}, along with an object-oriented interface
for specifying and solving resource allocation problems of the form
\eqref{e-res-alloc}. In our library, users can select from a library of utility
functions (or define their own). Our code is open-source, and available at

\begin{center}
\url{https://github.com/cvxgrp/resalloc}.
\end{center}

\paragraph{Solving the subproblems in parallel.} Our implementation of
algorithm~\ref{alg-max-net-utility} is completely vectorized and exploits the
fact that for each subproblem, there is a basic feasible solution in which at
most two resources are used. In particular, we compute the slopes of the affine
functions on which $c(t)$ depends using vectorized operations (across all
jobs), and collect them into a matrix of shape $n$ by $m(m-1)/2$ (this is
tractable, because in the problems we are interested in $m \ll n$). We then
operate on this matrix in a vectorized fashion to obtain the allocation solving
the subproblem.

Roughly speaking, this means we solve the $n$ subproblems in
parallel, but not by spawning multiple threads that solve the problems in isolation.
Instead, we make heavy use of vector and matrix operations (avoiding control
flow such as for loops as much as possible), and let our numerical linear
algebra software (\ie, PyTorch) exploit the parallelism that is intrinsic to
these operations. On a CPU, this means we exploit parallelism at multiple
levels: the vector and matrix operations are split over multiple
threads, and each thread in turn can take advantage of SIMD (single instruction,
multiple data) operations supported by the hardware
\cite[Chap.~4]{hennessy2011computer}. Additionally, our software library
supports CUDA acceleration via PyTorch, \ie, users can run the price discovery
method on GPUs; in this case, the vector and matrix operations in our
implementation are split over the GPU's streaming multiprocessors, each
of which can be thought of as a SIMD processor.

Because we efficiently exploit parallelism, our implementation is often orders
of magnitude faster than off-the-shelf solvers for convex optimization, and it
can scale to much larger problems. We will see some experiments that
demonstrate this in the following subsections.

\paragraph{Utility functions.} 
Our implementation is modular, and can support any utility function that
implements three specific methods: one that evaluates the utility, one
that solves the subproblem \eqref{e-opt-over-t} (given the slopes of the affine
functions on which $c(t)$ depends), and one that computes an initial guess for
the price vector (this guess need not be intelligent). These methods must be
vectorized over jobs, which is simple to do using PyTorch.

\paragraph{Code example.} Below we show a code example
that formulates a simple resource allocation problem using our software. (Here,
the entries of the resource limits and throughput matrix are randomly
generated; in an application, a user would use their real data for these
variables.)
\clearpage
\lstset{language=mypython}
\begin{lstlisting}
import torch
from resalloc.fungible import AllocationProblem, utilities

n_jobs, n_resources = int(1e6), 4
throughput_matrix = torch.rand((n_jobs, n_resources))
resource_limits = torch.rand(n_resources) * n_jobs + 1e3

problem = AllocationProblem(
  throughput_matrix=throughput_matrix,
  resource_limits=resource_limits,
  utility_function=utilities.Log()
)
\end{lstlisting}

The problem can then be solved by calling the solve method:

\lstset{language=mypython}
\begin{lstlisting}
problem.solve(verbose=True)
\end{lstlisting}
This yields the following verbose output, tracking the progress of the
algorithm; each line prints the average utility, dual function value (divided
by the number of jobs), and the gap between the two at a specific iteration
of the price discovery algorithm. (By default the solver uses L-BFGS for
the price updates, and terminates when the gap is less than $10^{-3}$.)
\lstset{language=mypython}
\begin{lstlisting}
iteration 00 | utility=-0.349803 | dual_value=0.275906 | gap=6.26e-01
iteration 05 | utility=-0.274916 | dual_value=-0.26548 | gap=9.44e-03
Converged in 009 iterations, with residual 0.000161
\end{lstlisting}
After calling the solve method, the optimal allocation can be obtained by
accessing the \texttt{X} attribute of the problem object (\texttt{problem.X}),
and the optimal prices can be obtained via the \texttt{prices} attribute
(\texttt{problem.prices}).

Solving this problem using a GPU requires just two changes to the code sample,
shown below.
\lstset{language=mypython}
\begin{lstlisting}
problem = AllocationProblem(
  throughput_matrix=throughput_matrix.cuda(),
  resource_limits=resource_limits.cuda(),
  utility_function=utilities.Log()
\end{lstlisting}

\paragraph{Experiment set-up.}
In the following subsections, we demonstrate our implementation on 
numerical examples. All experiments use our implementation with default
parameter values and 32-bit floating points.  We use L-BFGS with memory $10$
for the price updates, and we run the algorithms with
tolerance $\epsilon = 10^{-3}n$, which is more 
accuracy than needed in any practical application.
We run the algorithm on a CPU, an Intel i7-6700K CPU with four physical cores 
clocked at 4 GHz, and also
a GPU, an NVIDIA GeForce GTX 1070 with 1920 cores clocked at 1.5 GHz and a peak
of 6.5 TFLOPs. We also give some experiments on 
a more compute-intensive configuration, an NVIDIA DGX-1 equipped with an
NVIDIA Tesla V100 SXM2 GPU, which has 5120 cores clocked at 1.29 GHz and 
a peak of 15.7 TFLOPs.

Where possible, we cross check each solution found by our method 
using CVXPY \cite{diamond2016cvxpy, agrawal2018rewriting} with the MOSEK solver
\cite{mosek}, a high performance commercial interior-point solver. (CVXPY is
unable to compile one very large instance, and we suspect MOSEK would be unable
to solve it on our machine.) In all cases, the allocations found using our
method and MOSEK agreed.

\subsection{A medium size problem}
We consider a medium size problem, with $n=10^6$ jobs and $m=4$ resources.
The data are synthetic but realistic.
We choose the entries of $a_i$ from uniform distributions,
\[
  (a_i)_1 \sim U(0.1, 0.3) \quad (a_i)_2 \sim U(0.1, 0.5) \quad (a_i)_3 \sim U(0.3, 0.8) \quad (a_i)_4 \sim U(0.6, 1.0)
\]
and we choose $R$ as
\[
  R =  (8 \times 10^5, 10^5, 10^4 , 10^3),
\]
so that the resources that are more efficient (on average) are also more scarce.
We consider two problems, one with log utility and another with target-priority
utility, with $t_i = 0.2$ and priorities randomly chosen to be $w_i = 1$ or
$w_i = 2$ (each with probability one half).

We ran the price discovery algorithm on these problems, using the initialization
\eqref{e-initialization} for the prices. Maximizing the log utility took 3.9
seconds
on our CPU and 0.30 seconds on our GPU; maximizing the target-priority utility
took 11 seconds on our CPU and 0.89 seconds on our GPU.  These are impressive
solve times, considering our problem has 4 million variables.  
For comparison, MOSEK (which solves to high accuracy) took 300 seconds and 36
iterations for the log utility and 99 seconds and 74 iterations for the
target-priority. MOSEK's set-up time for the problems was 18 seconds and 13
seconds, respectively, and it it took MOSEK approximately 111 seconds and 85
seconds to reach solutions of the same accuracy as ours.

\begin{figure}
  \includegraphics{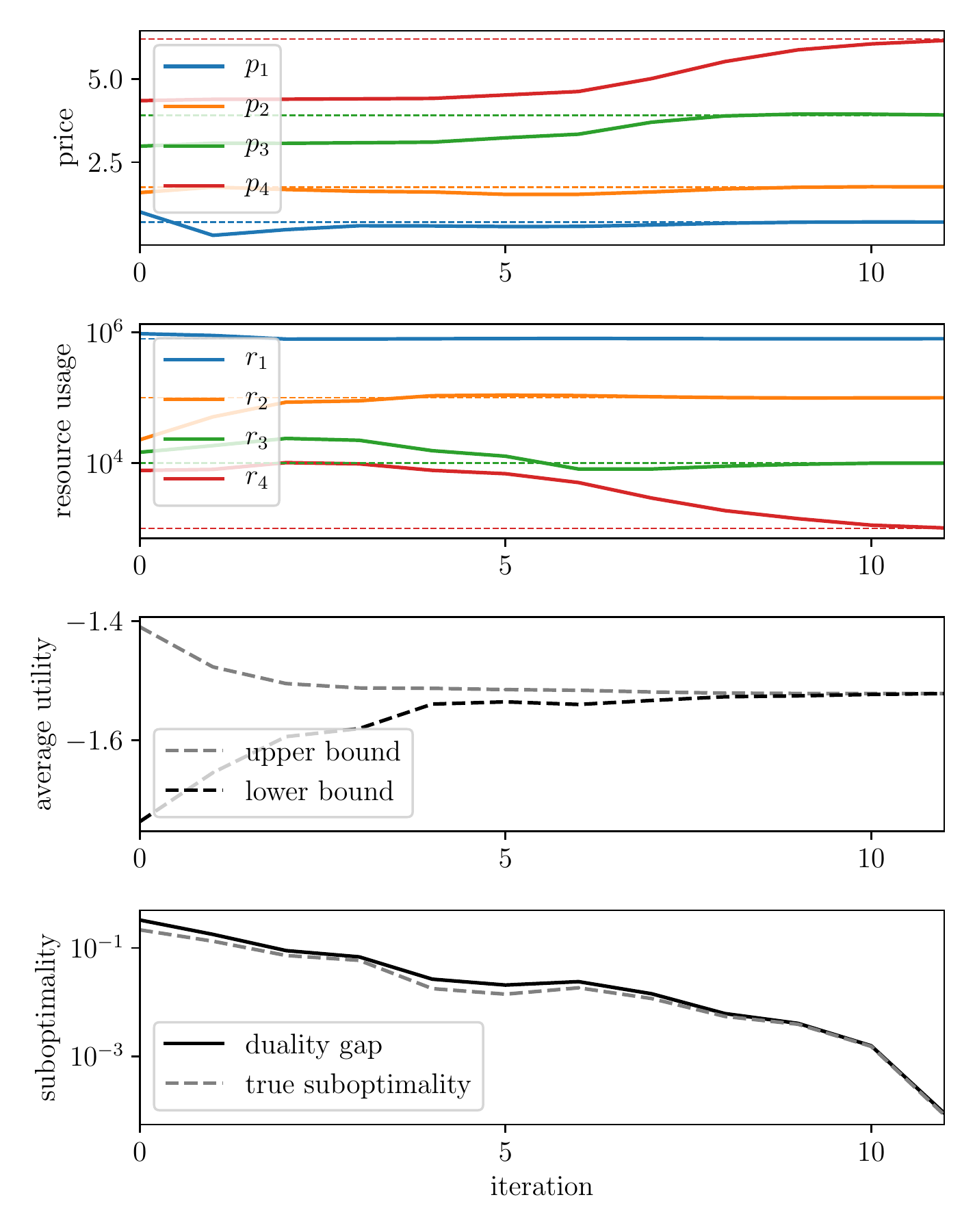} \caption{\emph{Log utility.} Prices,
resource usage, bounds on average log utility, and suboptimality
versus iteration number.}
\label{f-progress-log}
\end{figure}

\begin{figure}
  \includegraphics{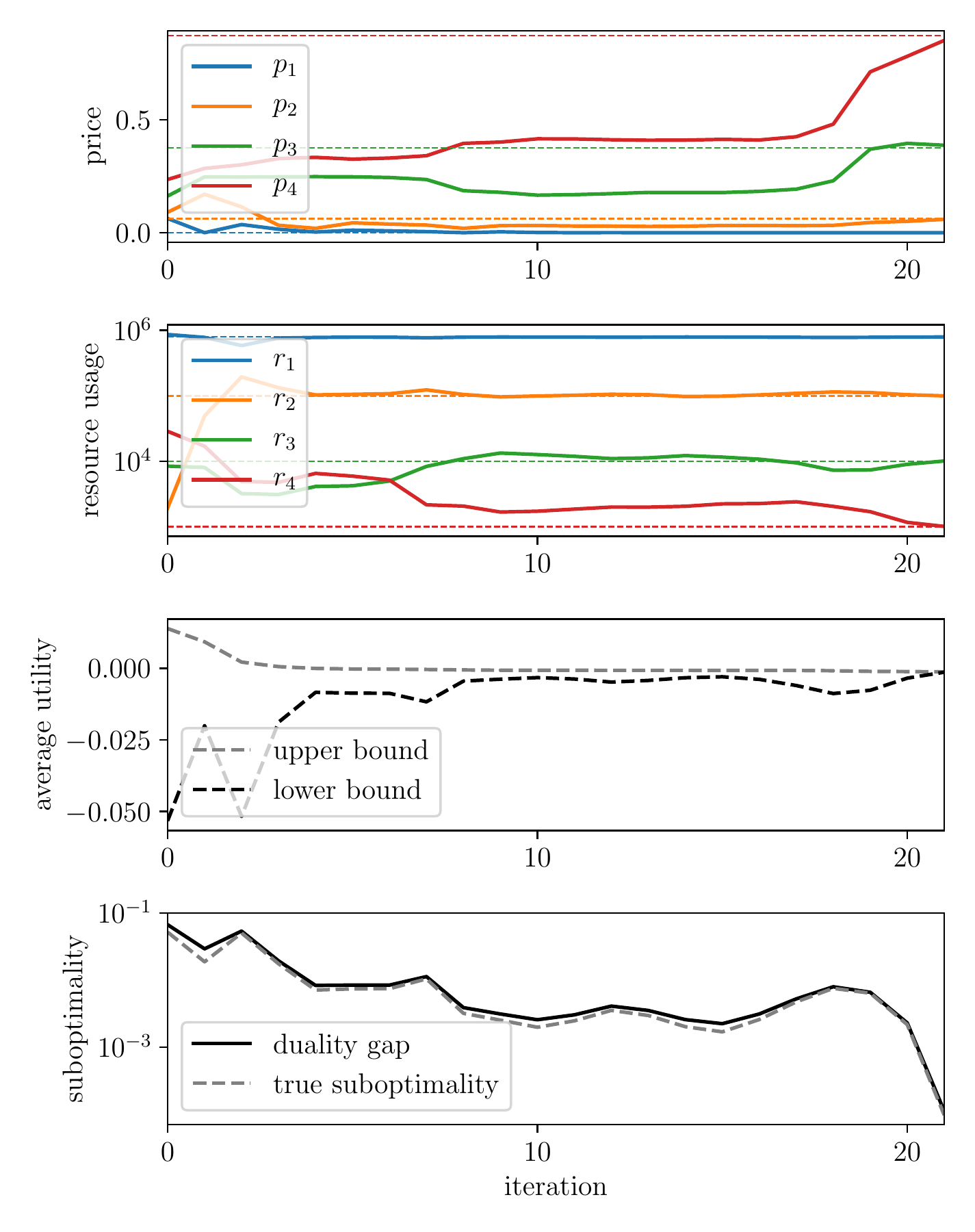} \caption{\emph{Target-priority
    utility.} Prices, resource usage, bounds on average target-priority
    utility, and suboptimality versus iteration number.}
\label{f-progress-target}

\end{figure}

The progress of the algorithm is shown in figures~\ref{f-progress-log} and
\ref{f-progress-target}, for the log and target-priority utilities,
respectively. In each figure the top plot shows the prices $p^k$ versus
iteration, with dashed lines showing the optimal prices $p^\star$.
We can see that the prices converge to optimal (within $10^{-3}$) in
11 iterations for log utility and 21 iterations for target-priority
utilities.
The second plot shows the resource usage $r^k$ versus iteration, with 
dashed lines showing the resource limits. The resource usage
constraints are violated early on but are satisfied by the end, as the prices
converge to their optimal values. The third plot
shows the upper and lower bounds on average
utility, and the bottom plot shows the duality gap divided by $n$,
$(g(X^k) - U(\tilde X^k))/n$ 
(which is known at iteration $k$) and the true suboptimality 
divided by $n$,
$(U^\star - U(\tilde X^k))/n$ (which is not
known at iteration $k$). For $U^\star$ we used the optimal value obtained by
MOSEK.

\begin{figure}
  \centering
\includegraphics{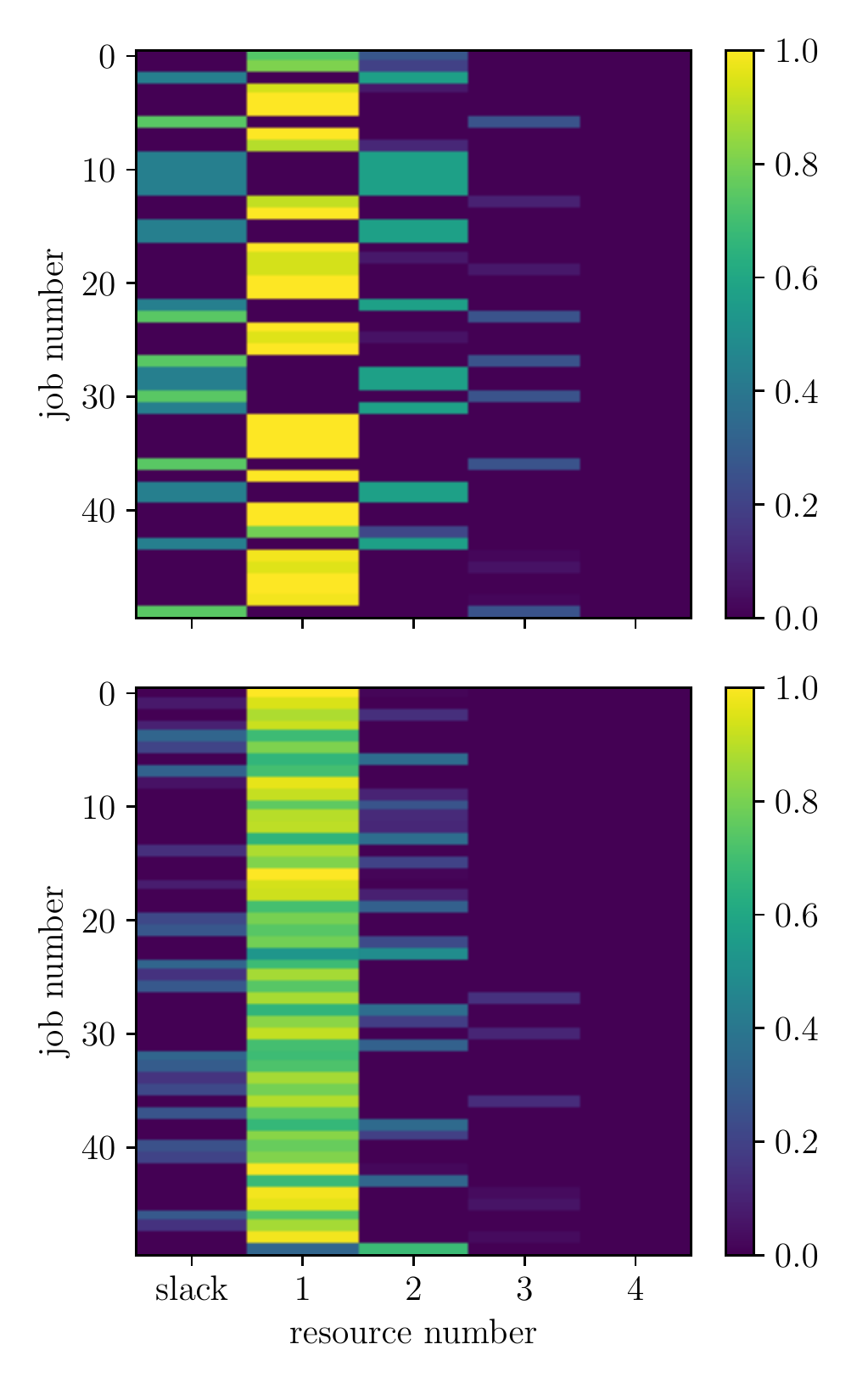}
\caption{Fifty rows of the optimal allocation matrices, along with the
  slack. Each row represents an allocation vector $x_i$, and each column
  a resource (or slack). Each job uses at most two resources. \emph{Top.}
  Log utility. \emph{Bottom.} Target-priority utility. }
\label{f-allocation-matrices}
\end{figure}

Figure~\ref{f-allocation-matrices} visualizes 50 rows of the optimal allocation
matrices, along with the slack; the rows were selected by embedding the
efficiency vectors $a_i$ into a line (using principal component analysis, via PyMDE
\cite{agrawal2021minimum}, a Python package for embedding), sorting the
resulting embedding, and permuting the rows of the allocation matrix to match
the sort order. For the allocation maximizing log utility, around 17 percent of
jobs use two resources, the remaining use just one resource, and roughly 18
percent of jobs have a positive slack (meaning that $\ones^Tx_i <1$, \ie, these
jobs would run for less than 100\% total time). For the allocation maximizing
the target-priority utility, 67 percent of jobs use two resources, the
remaining use just one, and 50 percent have
positive slack.  

\begin{figure}
\includegraphics{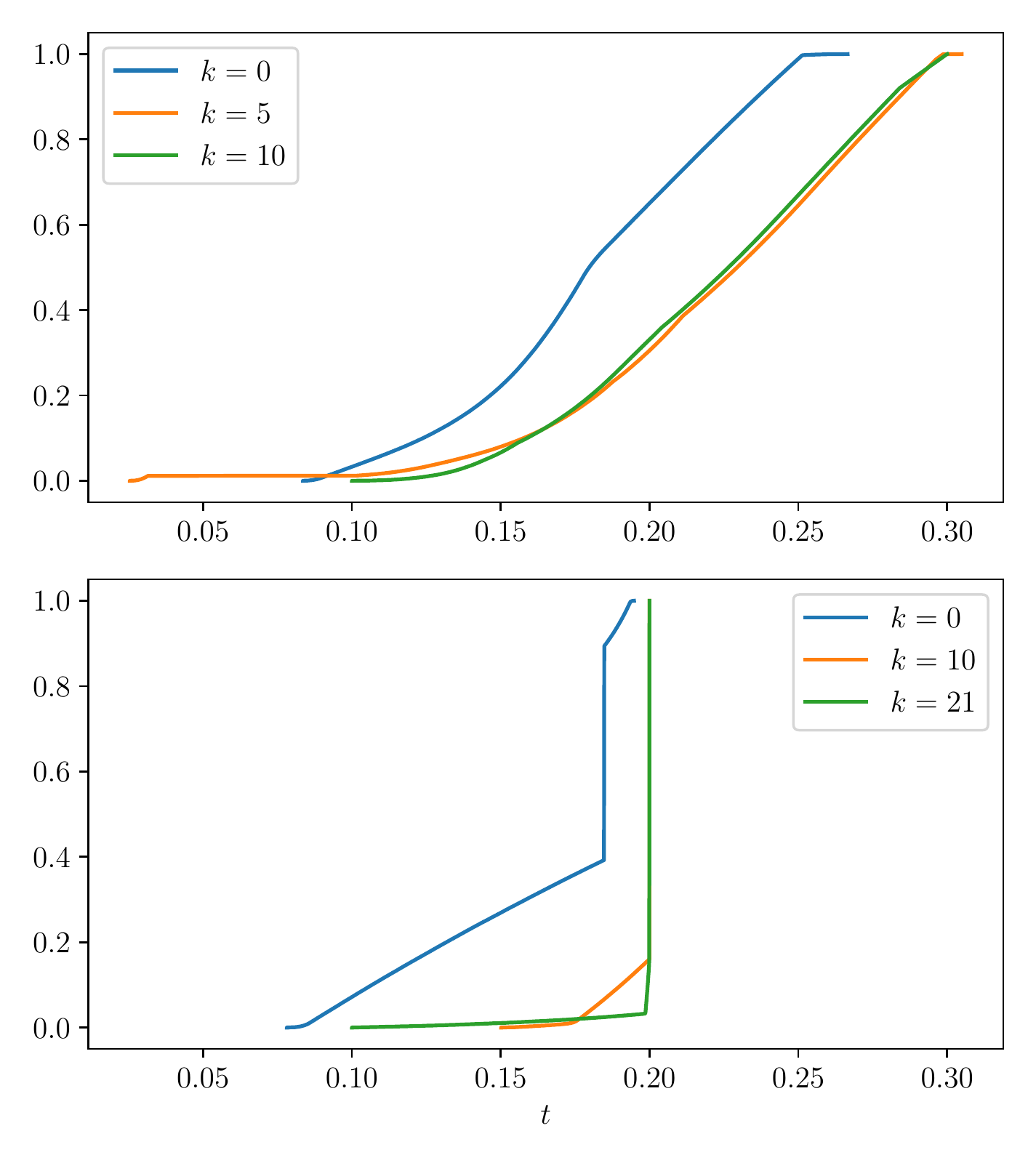}
\caption{CDFs of throughputs, at various iteration numbers $k$.
\emph{Top.} Log utility. \emph{Bottom.} Target-priority utility.}
\label{f-throughput-distributions}
\end{figure}

\begin{figure}
\includegraphics{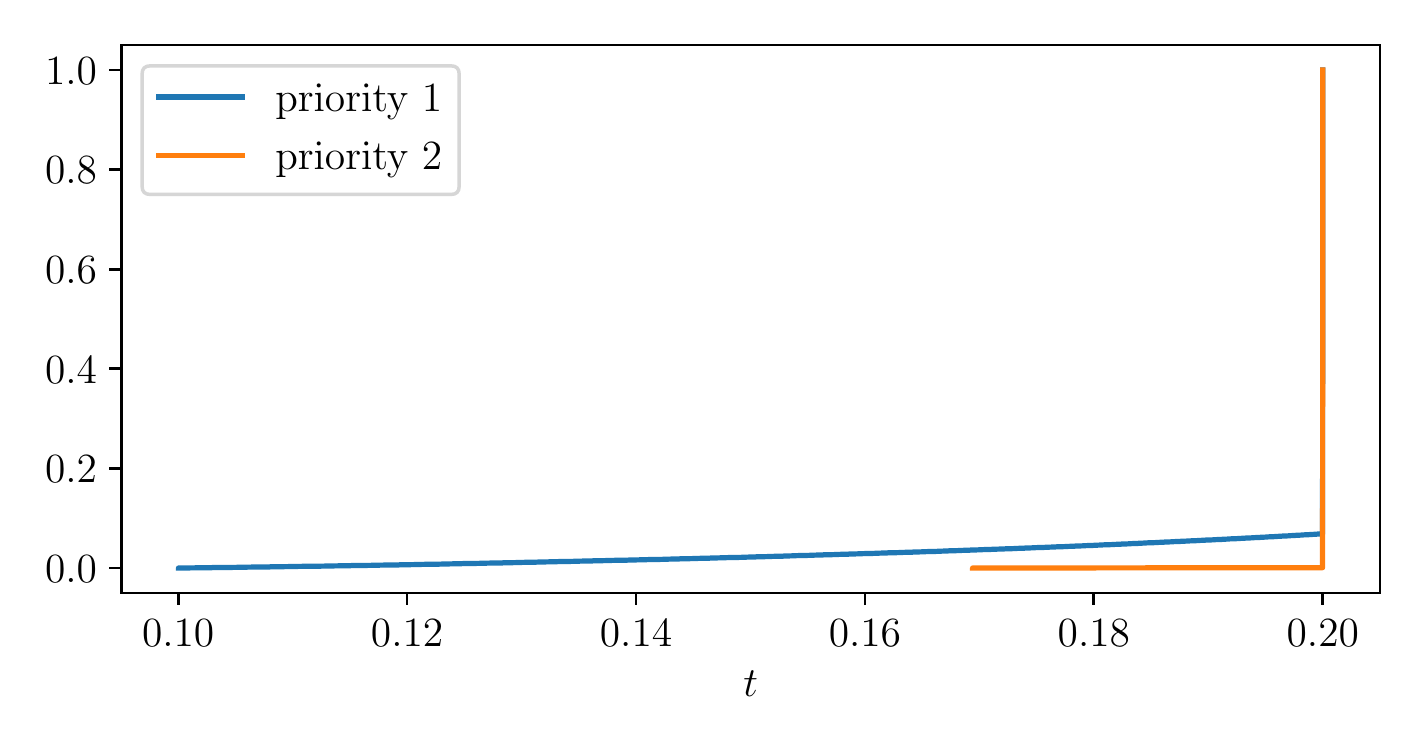}
\caption{CDFs of throughputs for optimal target-priority allocations,
for the low and high priority jobs.}
\label{f-throughput-distribution-vs-priority}
\end{figure}

The job throughput
distributions for the initial allocation, the allocation 
obtained after a few iterations,
and the final optimal allocation
$X^\star$ are shown for the two problems 
in figure~\ref{f-throughput-distributions}.
We can see that the target-priority allocation is able to obtain the target
throughput $0.20$ for over 95\% of the jobs.
Figure~\ref{f-throughput-distribution-vs-priority} shows
the throughput distributions for the target-priority optimal allocation,
for the jobs with low priority ($w_i=1$) and higher priority ($w_i=2$),
respectively.
Virtually all of the high priority jobs achieve the target throughput,
while around 10\% of the low priority jobs do not.
(Recall that roughly half the jobs are high priority, and half
are low priority.)


\clearpage

\subsection{A large problem}
For our next example we consider a large problem with $n=50$ million jobs,
$m=4$ resources, and log utility. This is a
convex optimization problem with 200 million variables. We generate the
entries of $a_i$ in the same way as the previous example, and use the previous
resource limits scaled by 50.

\begin{figure}
\includegraphics{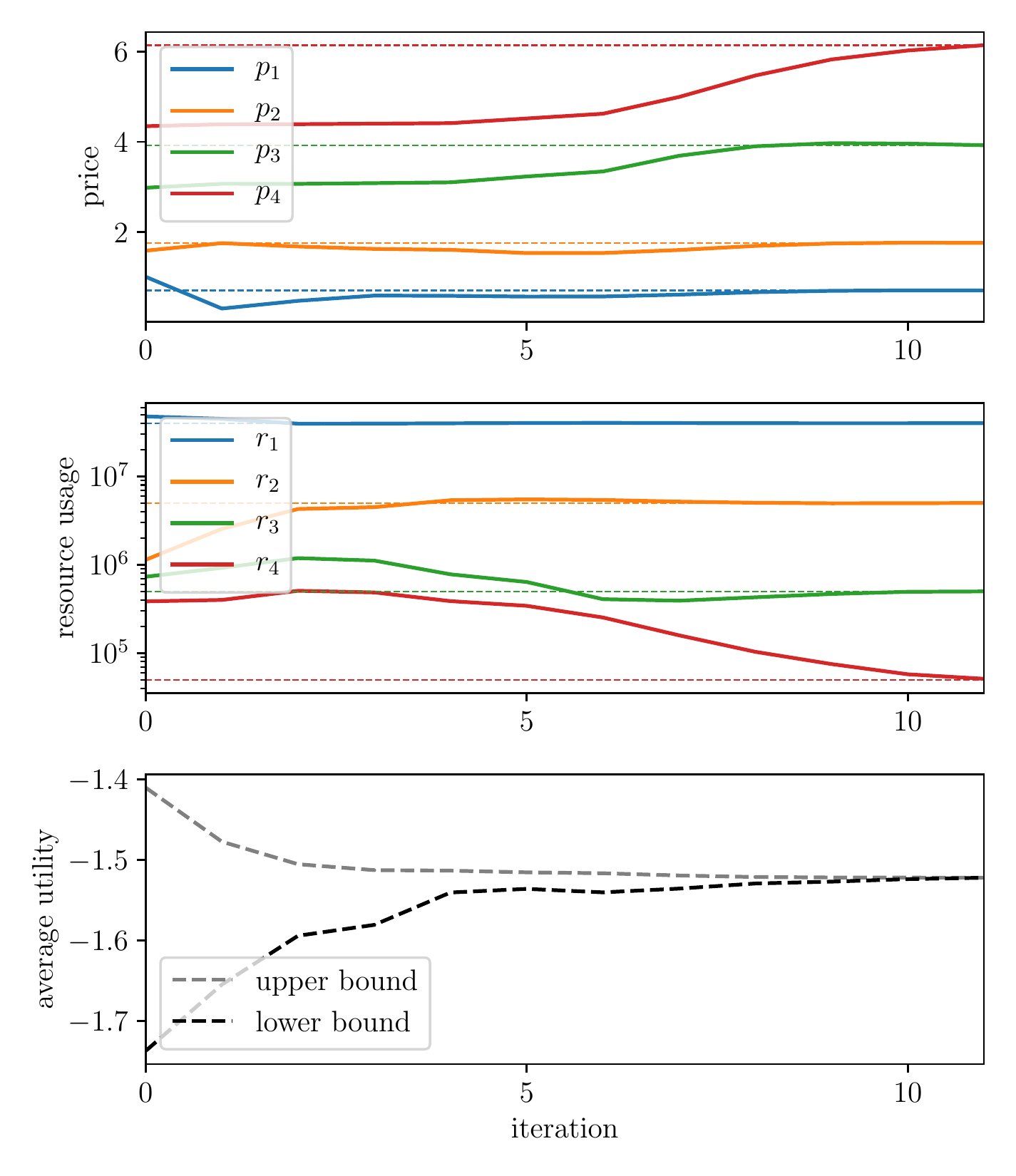} \caption{\emph{A large
problem.} Prices, resource usage, and bounds on average log utility
versus iteration number, for a large problem with $n=50$ million jobs and $m=4$
resources. This problem was solved in 228 seconds on  CPU and 7.9 seconds on
GPU.} \label{f-progress-log-big}
\end{figure}

We solved this problem on our smaller machine's CPU, and on the Tesla V100 GPU
(the problem was too large to fit on the smaller GPU). The problem was solved
in 228 seconds on our CPU, and in roughly 7.9 seconds on the GPU. (This problem
was too large to solve with MOSEK.)  Figure~\ref{f-progress-log-big} shows the
progress of the algorithm. In the final allocation, 17 percent of jobs use two
resources, the remaining use one, and 18 percent have positive slack.

\clearpage

\subsection{Scaling with jobs and resources}

In this subsection we show how our method scales in the number of jobs and
resources, and compare solve times for a number of problems against
low-accuracy MOSEK solves.

\begin{figure}
\includegraphics{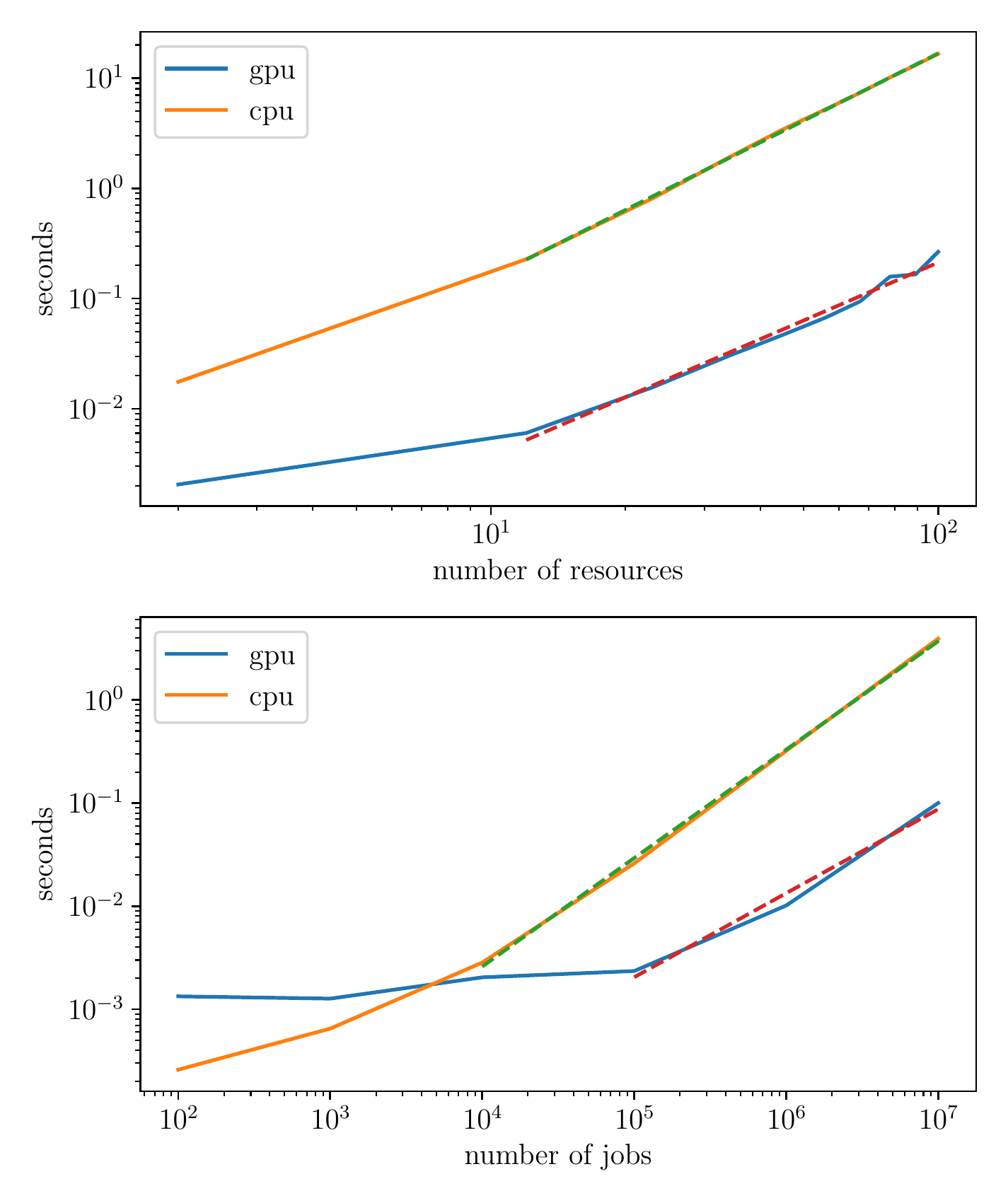}
\caption{\emph{Evaluating the dual function.} Mean time elapsed evaluating the
  dual function, with dashed lines showing scaling over
  ranges where the fixed overhead no longer dominates. \emph{Top.} $n=10^5$,
  with $m$ varying. \emph{Bottom.} $m=4$, with $n$ varying.}
\label{f-scaling}
\end{figure}

\paragraph{Evaluating the dual function.}
The main computation in each iteration of our algorithm is evaluating
the dual function (\ie, solving the $n$ subproblems). To give an idea of how our
implementation scales, we timed how long it took to evaluate the dual function
for synthetic problems, varying the number of jobs $n$ and resources $m$. We
conducted two experiments. In one, we held $n$ fixed at $10^5$, and varied $m$
from $2$ to $100$; in the other, we held $m$ fixed at $4$, and varied $n$ from
$10^2$ to $10^7$. For each instance we evaluated the dual function five times.
The mean times are plotted in figure~\ref{f-scaling}. To study the
scaling, for each plot we log-transformed the inputs and outputs and fit
linear regressions, \ie, we fit models of the form
\[
  \log(s) \approx a + b \log(m), \quad \log(s) \approx a + b \log(n),
\]
where $s$ is the elapsed time. We fit these models in ranges where the
constant overhead of our software implementation no longer dominated.
These fits are plotted as dashed lines in figure~\ref{f-scaling}.

For CPU, the coefficient $b$ is 2.03 for resources and 1.05 for jobs; for GPU
(we used the Tesla V100), $b$ is 1.75 for resources and 0.81 for jobs. (The
coefficients $a$ are very small: for CPU, $\exp(a)$ is $1.4 \times 10^{-3}$ and
$1.6 \times 10^{-7}$ for resources and jobs respectively, and for GPU the
values are $6.7 \times 10^{-5}$ and $1.7 \times 10^{-7}$.) This means that CPU
time scales quadratically in the number of resources and linearly in the number
of jobs, while GPU times scales less than quadratically in the number of
resources and sublinearly in the number of jobs, at least over the range
considered.

\clearpage
\paragraph{Comparison to MOSEK.}
For large problems, our method outperforms high-quality off-the-shelf solvers
for convex optimization, such as MOSEK. To demonstrate this, we timed how long
it took our method to solve a number of synthetic problems (with log utility
and accuracy $\epsilon=10^{-3}n$), varying the number of jobs $n$ and resources
$m$; we solved the same problems with MOSEK, to low accuracy (setting all
tolerances to $10^{-3}$). In the first set of experiments, we held $n$ fixed at
$10^6$, and varied $m$ from $2$ to $16$; in the second, we held $m$ fixed at
$4$, and varied $n$ from $10^2$ to $10^6$. The columns of the throughput
matrices were sampled from uniform distributions, so that (on average) the
resources were ordered from least to most efficient. The resource limits were
generated by sampling each entry from a uniform distribution on $[0.1, 1]$, and
scaling the first component by the number of jobs, the second by the number of
jobs divided by 1.5, the third by the number of jobs divided by $1.5^2$, and so
on. We solved five instances of each problem.

\begin{figure}
\includegraphics{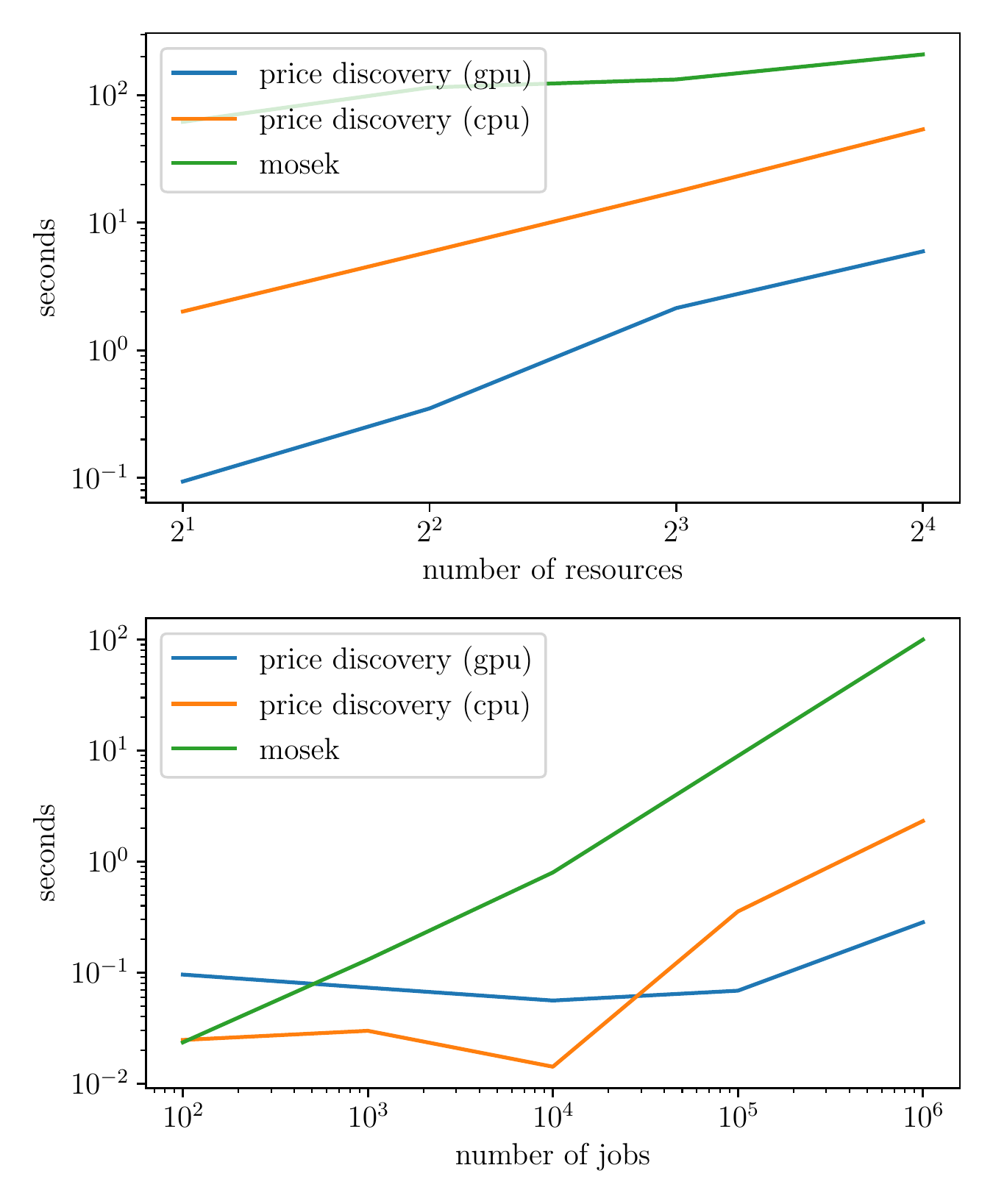}
\caption{\emph{Comparison to MOSEK.} Mean time elapsed solving resource
  allocation problems with log utility, comparing our method and low-accuracy
  MOSEK solves. \emph{Top.} $n=10^6$, with $m$ varying. \emph{Bottom.} $m=4$,
  with $n$ varying.}
\label{f-v-mosek}
\end{figure}

The mean solve times are plotted in figure~\ref{f-v-mosek}. For 
large problems, our method appears to be between one to three orders of
magnitude faster than low-accuracy MOSEK solves. We point out that we solve the
largest instance, which has 16 million variables ($n=10^6$ jobs and $m = 16$
resources), in just 5 seconds on a GPU. On this same instance, MOSEK takes over
200 seconds.

\clearpage

\section{Conclusion}\label{s-conclusions}

We have described a custom solver for the fungible resource allocation
problem that scales to extremely large problem instances, especially 
when run on a GPU.  For example, problems with millions of variables
can be solved in well under a second.  Smaller problems can be solved in
milliseconds.

The method uses the dual problem, and manipulates a set of resource prices
until we achieve an optimal allocation.  The optimal prices, which our 
algorithm discovers, can be used for other tasks as well as a solving
the allocation problem.  For example we can actually charge jobs in proportion
to $(p^\star)^T x_i$, the total cost of the resources used with the optimal
prices.

The prices can also be used in situations where an allocation problem is broken
into $M$ shards or subproblems, and each of these solved independently with
resource budget $(1/M)R$ (as proposed in \cite{narayanan2021dont}). If the
optimal resource prices for the different shards are close, we can conclude
that the solution found from the partitioned problems is nearly optimal for the
allocation problem, had they been solved together. If the resource prices vary
across the shards, they can be used to re-allocate resources to the shards, by
moving resources from the shards with lower prices to those with higher prices.
This method, which is very interpretable, can be shown to converge to the
solution of the larger problem \cite[\S3.1]{boyd2007notes}, and it can be used
to improve suboptimal solutions obtained via the method from
\cite{narayanan2021dont}. This method could also be used to re-allocate
resources across separate virtual data centers in a principled way.

\subsection{Extensions and variations}
We mention a few extensions and variations on our method, and related
problems. The first is the min-throughput utility. This objective function is
concave and 
nondecreasing, but it is not separable, so the methods of this paper cannot
be directly used.  But very similar methods can be. Methods described in
this paper can be used to approximate the min-throughput utility; for example,
just using a utility function with strong curvature, like a power utility with
$p$ near zero, already gives a good approximation of the minimum throughput.
Or we can use a target-priority utility and decrease the target throughput
until all (or a large fraction) of the jobs meet the target.

Our formulation \eqref{e-res-alloc} can easily be extended to the case
in the constraint on the resource usage is replaced with
$\sum_{i=1}^n d_i x_i \leq R$, where $d_i \in \reals_{++}$ is the number
of resources demanded by job $i$ (this constraint was proposed in
\cite{narayanan2020heterogeneity}). For example, if $d_i = 2$, and the
resources are GPUs, this means that job $i$ requires two GPUs in order to run.
We can also handle the case in which jobs demand different amounts of
different resources (leading to demands that are vectors in $\reals^m_{++}$, not
scalars). These extensions would require very minimal changes to our price
discovery method (the only change is that the price vector in each subproblem
\eqref{e-subproblem} is multiplied by the job's demand); indeed, we have
implemented them in our software.

In this paper we only considered the problem of allocating the fraction
of time that a job should spend consuming each resource; we did not consider
the problem of coming up with a schedule, which says when each job should
consume its resources during the time interval. In our
motivating application of computer systems, some pairs of jobs
exhibit poor performance when colocated on the same hardware (\ie, when
using the same resource) \cite{kambadur2012measuring, zhang2013cpi2, verma2015large,
narayanan2020heterogeneity}. This kind of interference could be mitigated by a
lower-level scheduler, which takes as input an optimal allocation matrix, and
comes up with a concrete schedule that respects the time-slicing while
minimizing cross-job interference. The lower-level scheduler could
itself be based on solving an optimization problem; Netflix's optimization-based
scheduler is one such example \cite{rostykus2019predictive}. We mention
that it is also possible to take interference into account in a higher-level
optimization-based scheduler, possibly using a formulation of
the resource allocation problem that is different from ours (one such
formulation is given in \cite[\S3.1]{narayanan2020heterogeneity}).

Finally, another related and important problem is the allocation of
non-fungible resources. This occurs when the resources are heterogeneous, and
not fully exchangeable. For example in a data center, we can allocate number of
cores, I/O bandwidth, memory, and disk space. (These are not fungible, because
you cannot get any throughput if you only use I/O bandwidth and no cores,
memory, or disk space.) In this setting, the utility function for each job does
not have the form of a scalar utility function of a linear function of the
resources allocated; instead it is a utility function of $x_i$. These
problems can often be posed as convex optimization problems (\eg,
\cite{ghodsi2011dominant, bird2011pacora}.) The specific solution of the
subproblem described in \S\ref{s-subprob} no longer applies, but the price
discovery method in general does. We will address that problem in a future
paper.

\clearpage

\printbibliography
\end{document}
\printbibliography
\end{document}